\documentclass[10pt,twocolumn,letterpaper]{article}

\usepackage{cvpr}              

\definecolor{cvprblue}{rgb}{0.21,0.49,0.74}
\usepackage[pagebackref,breaklinks,colorlinks,allcolors=cvprblue]{hyperref}
\usepackage[ruled,vlined]{algorithm2e}

\usepackage{float}       
\usepackage{placeins}    
\usepackage{graphicx}    
\usepackage{caption}     
\usepackage{float}
\usepackage{tabularx}
\usepackage{array}
\usepackage{enumitem}

\def\paperID{} 
\def\confName{}
\def\confYear{2025}

\title{Understanding Student Perceptions of Flipped Linear Algebra Classrooms\\ via Interpretable Machine Learning}

\author{
Sudip Laudari\\
{\normalsize Department of Civil Engineering} \\
{\normalsize The University of Sydney, New South Wales 2006, Australia} \\
{\tt\small sudip.laudari@sydney.edu.au}
\and
N. Karjanto\\
{\normalsize Department of Mathematics, College of Natural Science} \\
{\normalsize Sungkyunkwan University, Suwon 16419, Gyeonggi-do, Republic of Korea} \vspace{0.1cm}\\
{\normalsize Humanities, Arts, and Social Sciences Division, Underwood International College}\\
{\normalsize Yonsei University, Songdo, Incheon 21983, Republic of Korea} \\
{\tt\small natanael@skku.edu}
}

\begin{document}
\maketitle
\begin{abstract}
Flipped classroom pedagogy is widely used in undergraduate mathematics to promote active learning, yet it remains unclear whether students experience it in systematically different ways. In this study, we analyze student perceptions from an introductory linear algebra course using survey data collected across multiple semesters. Using an interpretable machine learning approach, we examine patterns across questionnaire responses and evaluate their consistency under repeated analysis. Our results reveal a clear and stable separation in perception patterns when grouped by gender, suggesting that these differences arise from structured combinations of factors rather than isolated responses. The model also identifies key aspects of engagement and instructional design that contribute most to this separation. These findings highlight opportunities for more inclusive flipped classroom design and demonstrate the value of interpretable methods in educational research.

\end{abstract}
    
\section{Introduction}
\label{sec:intro}

As educators, we aim to support students’ learning by promoting knowledge retention, problem-solving ability, and critical thinking skills. Active learning approaches have received increasing attention as alternatives to traditional lecture-based instruction. The flipped classroom has emerged as one of the most widely adopted pedagogical models within this framework \cite{bergmann2023flip, bishop2013flipped, lage2000inverting}. In a flipped classroom, students review new material before class through videos or readings. This allows in-class time to be devoted to discussion, problem solving, and other higher-order learning activities.

Active learning strategies have been shown to enhance instructional effectiveness and student engagement across many disciplines. In mathematics education, and particularly in linear algebra, flipped classrooms have been implemented using a variety of instructional designs \cite{talbert2014inverting}. Empirical comparisons between flipped and traditional instruction frequently report improvements in conceptual understanding and learning attitudes \cite{love2014student, murphy2016student}. Similar findings have been observed across different institutional contexts, including increased participation and improved self-directed learning skills \cite{novak2017flip, nasir2020effectiveness, karjanto2019english, karjanto2022sustainable}. These outcomes are also consistent with constructivist theories which emphasize active knowledge construction and peer interaction \cite{prince2004does}.

Despite strong evidence supporting flipped classrooms, much of the existing literature focuses on aggregate outcomes such as average examination performance \cite{strelan2020flipped,he2016effects}, overall course achievement \cite{love2014student,murphy2016student,karjanto2022sustainable}, and broad measures of student satisfaction \cite{strelan2020student,o2015use}. Such analyses often treat student populations as homogeneous groups and may obscure differences in how students experience the same instructional design. Relatively little attention has been paid to whether different groups of students experience flipped learning in systematically different ways. Differences related to gender or prior preparation may remain hidden when analyses rely primarily on aggregate statistics. Understanding these differences is important for refining pedagogical design and ensuring that active learning remains inclusive \cite{freeman2014active,eddy2014getting}.

From a methodological perspective, survey-based studies in education are commonly analyzed using item-by-item comparisons. While these approaches provide useful summaries, they may overlook multivariate relationships in which perceptions are shaped by combinations of factors rather than isolated responses \cite{whitelock2019student,ranjeeth2020survey}. Machine-learning methods provide a complementary framework by enabling the identification of structured patterns in high-dimensional survey data \cite{romero2020educational}. These tools support a more nuanced, data-driven analysis of subgroup differences.

Accordingly, we adopt an interpretable machine-learning approach based on a linear support vector machine (SVM) to examine student responses \cite{cortes1995support,noble2006support,suthaharan2016support}. Our goal is to assess whether student perception patterns exhibit a stable and interpretable separation when stratified by gender. To support this exploratory analysis, we also employ principal component analysis (PCA) \cite{jolliffe2011principal}. The use of a linear model allows the results to be directly linked to individual questionnaire items for pedagogical interpretation \cite{livieris2023advanced}.

In this study, flipped classroom pedagogy is implemented in an introductory linear algebra course at a private university in South Korea. Student perceptions are collected through an end-of-semester questionnaire consisting of ten Likert-scale items. This work investigates whether interpretable machine-learning methods can reveal meaningful differences in how flipped environments are perceived across student groups. These insights aim to support more inclusive and effective pedagogical design.

\section{Methodology}
\label{sec:methodology}

\subsection{Data Collection}

Student perception data were collected using a structured questionnaire administered at the end of an introductory linear algebra course delivered using a flipped classroom format. Conducting the survey at the end of the semester allowed students to reflect on their overall learning experience after completing all instructional components, including pre-class materials, in-class activities, and assessments.

The questionnaire consisted of ten statements designed to capture students’ perceptions of flipped classroom pedagogy and related instructional practices. Responses were recorded on a five-point Likert scale:

\begin{itemize}[leftmargin=*]
\item \textbf{1 -- Strongly disagree}
\item \textbf{2 -- Disagree}
\item \textbf{3 -- Neutral}
\item \textbf{4 -- Agree}
\item \textbf{5 -- Strongly agree}
\end{itemize}

Students also provided basic demographic information including gender, year of study, academic major, and expected course grade. The original dataset contains approximately 200 student responses; however, the cohort exhibits an imbalanced gender distribution typical of STEM fields, with approximately 70--80\% male and 20--30\% female students.

Most flipped-classroom studies rely exclusively on authentic institutional data, such as examination results, survey responses, and learning management system (LMS) activity logs 
\cite{love2014student,murphy2016student, hardebolle2022gender,takrouni2024data}. Consistent with this practice, our analysis begins with real survey responses collected from enrolled students. However, educational cohorts are often relatively small, and class distributions may be imbalanced. To reduce sensitivity to limited sample size and class imbalance, we employ a model-based dataset regeneration strategy that preserves the empirical response patterns of the original cohort. The use of synthetic data to address small or imbalanced educational datasets has been explored in prior studies \cite{flood2025predicting,zhao2025large,santana2025exploring}.

To generate additional samples, response probabilities for each questionnaire item on the five-point Likert scale are estimated separately for male and female students using the original dataset. Synthetic samples are then created by first sampling a gender label according to the observed class proportions and subsequently sampling each of the ten questionnaire responses from the corresponding gender-conditional Likert distributions. This procedure preserves the observed response tendencies while increasing the effective sample size for analysis.

The regenerated samples are used only to stabilize model training and enable repeated stratified evaluation. They are not treated as independent empirical observations. The resulting dataset contains approximately 5{,}100 observations.

\begin{figure}[!t]
    \centering
    \includegraphics[width=\linewidth]{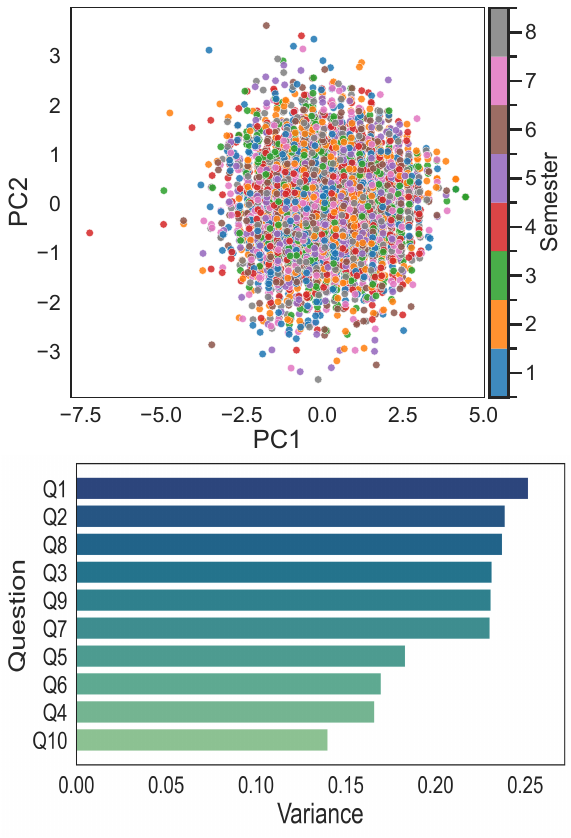}
    \caption{Principal component analysis of survey responses. 
    \textbf{Top:} Projection onto the first two components shows overlapping distributions. 
    \textbf{Bottom:} Variance across questionnaire items is uneven.}
    \label{fig:fig1}
\end{figure}

\subsection{Feature Encoding and Representation}

To prepare the questionnaire responses for linear modeling, each Likert-scale item is transformed into a centered numerical feature. Let $y_{ij} \in \{1,2,3,4,5\}$ denote the response of student $i$, for $i=1,\ldots,N$, to questionnaire item $j$, for $j=1,\ldots,10$. Each response is encoded as

\begin{equation}
x_{ij} = \frac{y_{ij} - 3}{2},
\end{equation}

which maps the original Likert scale to the symmetric interval $[-1,1]$. This transformation centers responses at the neutral midpoint and places all questionnaire items on a common numerical scale. Negative values represent disagreement, positive values represent agreement, and the magnitude reflects response strength. The encoding preserves the ordinal structure of the Likert scale while allowing all items to contribute comparably to the analysis.

The approach follows the common assumption that Likert responses approximate discretized observations of an underlying continuous perception variable with approximately equal spacing between adjacent categories \cite{carifio2008resolving,norman2010likert}. Under this representation, each student is described by a feature vector

\[
\mathbf{x}_i = (x_{i1}, x_{i2}, \ldots, x_{i10}) \in [-1,1]^{10}.
\]

A neutral response pattern corresponds to the origin of the feature space, providing a natural reference point for interpreting linear decision boundaries. Unlike sample-dependent standardization methods such as $z$-score normalization, this encoding preserves a fixed semantic baseline and allows direct interpretation of model coefficients as deviations from neutrality. 

\begin{figure*}[!t]
\centering
\includegraphics[width=\textwidth]{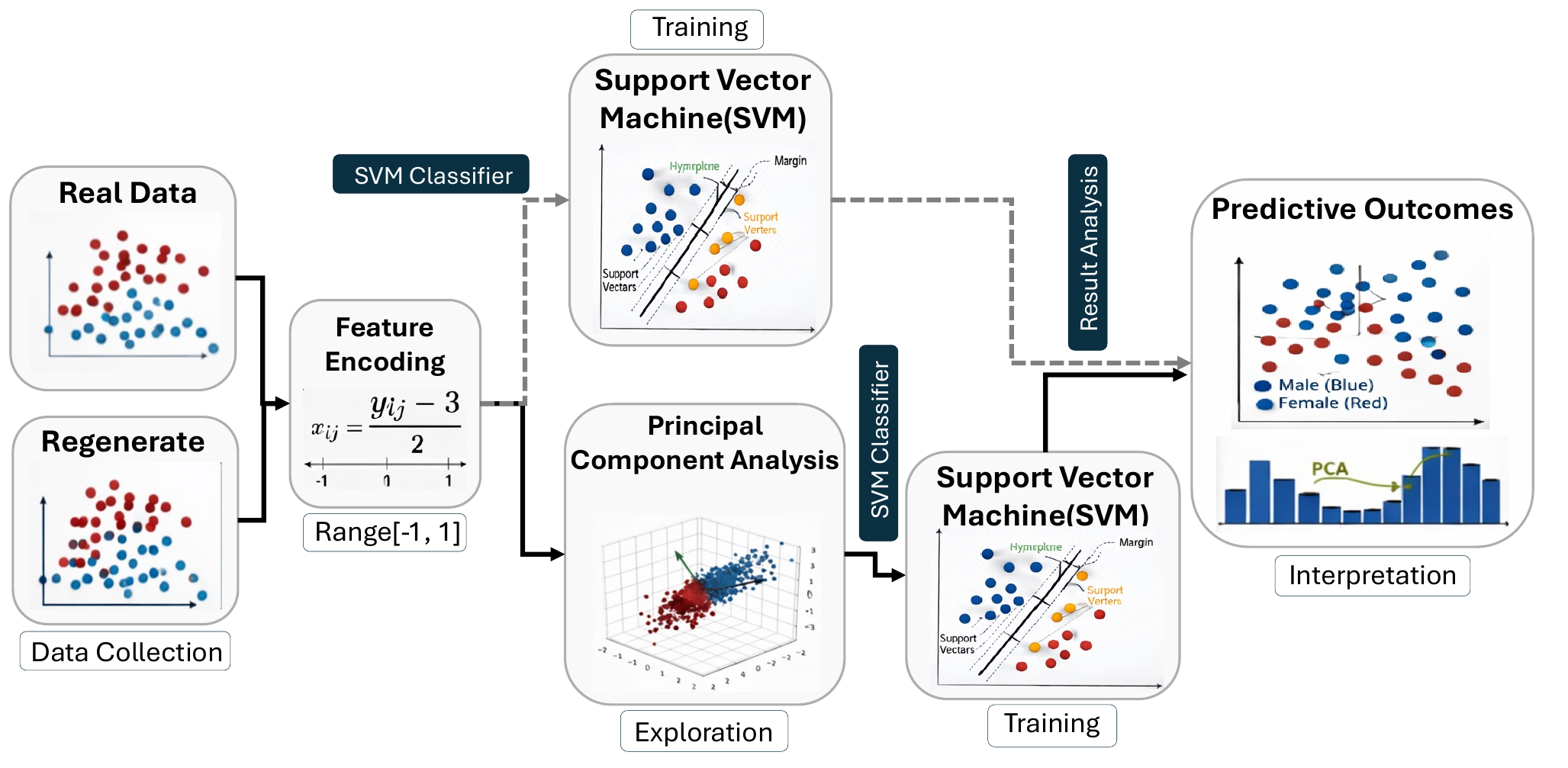}
\caption{Overview of the proposed framework. 
Survey responses are encoded to a centered scale and augmented for stability, then analyzed using PCA for exploratory structure and a linear SVM for classification, enabling interpretable identification of patterns in student perceptions.}
\label{fig:fig2}
\end{figure*}

The complete analysis workflow is summarized in Figure~\ref{fig:fig2}. Survey responses are first encoded into a centered numerical scale, regenerated samples are used to stabilize repeated evaluation, and the resulting feature representations are analyzed using PCA and a linear classifier.

\subsection{Exploratory Data Analysis via PCA}

Before classification, exploratory analysis is conducted to examine the structure of the encoded questionnaire data. PCA is applied to the centered feature vectors to identify dominant directions of variability and to support low-dimensional visualization. This analysis provides an interpretable projection of the data and enables comparison between models trained in the original feature space and those trained in a reduced representation.

Figure~\ref{fig:fig1} shows the projection of student responses onto the first two principal components and the variance contribution of individual questionnaire items. The projection exhibits substantial overlap among samples, indicating limited low-dimensional separation. The variance contribution plot shows that response variability is unevenly distributed across questionnaire items.


\subsection{Implementation of Machine Learning}

To analyze patterns in student perceptions, we employ a supervised classification framework based on a linear support vector classifier (SVC). The objective is not solely to maximize predictive accuracy, but to examine whether students’ response patterns exhibit a stable and interpretable linear separation when stratified by gender.

Let $\{(\mathbf{x}_i, y_i)\}_{i=1}^{N}$ denote the dataset, where $i=1,\ldots,N$, $\mathbf{x}_i \in \mathbb{R}^{10}$ represents the encoded questionnaire responses of student $i$, and $y_i \in \{0,1\}$ denotes the corresponding gender label. The linear SVM solves the following convex optimization problem:

\begin{equation}
\begin{aligned}
\min_{\boldsymbol{\omega}, b, \boldsymbol{\epsilon}}
&\quad
\frac{1}{2}\|\boldsymbol{\omega}\|^2 + C \sum_{i=1}^{N} \epsilon_i \\
\text{subject to}
&\quad
y_i (\boldsymbol{\omega}^{\top}\mathbf{x}_i + b) \geq 1 - \epsilon_i,
\quad i=1,\ldots,N, \\
&\quad
\epsilon_i \geq 0,
\quad i=1,\ldots,N.
\end{aligned}
\end{equation}

Here, $\boldsymbol{\omega}$ defines the normal vector of the separating hyperplane, $b$ is the bias term, and $\epsilon_i$ are slack variables that allow margin violations. The regularization parameter $C \in \mathbb{R}_{>0}$ controls the trade-off between maximizing the margin and penalizing classification errors.

The resulting maximum-margin solution defines a separating hyperplane in the feature space. Each component of the learned weight vector $\boldsymbol{\omega}$ corresponds to a specific questionnaire item, allowing direct interpretation of how individual survey statements contribute to the observed separation between gender groups.

To examine the influence of dominant variance structure on classification performance, models are trained both in the original feature space and in a reduced representation obtained through PCA. This comparison allows us to assess whether discriminative information aligns with the principal directions of variance identified during exploratory analysis.

\section{Training and Evaluation Protocol}
\label{subsec:training_protocol}

Each student response is represented as a 10-dimensional feature vector derived from the encoded questionnaire items (Q1--Q10). Model evaluation follows standard supervised learning practice using stratified train--test splits, cross-validation for hyperparameter selection, and held-out test evaluation \cite{kohavi1995study,pedregosa2011scikit}. Specifically, $70\%$ of the data is used for training and $30\%$ is reserved for testing. Stratification ensures that the gender distribution is preserved in both subsets, and a fixed base random seed is used to support reproducibility. The main training and evaluation settings are summarized in Table~\ref{tab:training_setup}.

\begin{table}[h]
\centering
\caption{Summary of model training and evaluation configuration.}
\label{tab:training_setup}
\renewcommand{\arraystretch}{1.1}
\begin{tabular}{@{}p{0.38\columnwidth} p{0.58\columnwidth}@{}}
\toprule
\textbf{Component} & \textbf{Description} \\
\midrule
Original dataset & $\sim$200 real student responses \\
Regenerated dataset & $\sim$5,100 observations (real + regenerated samples) \\
Number of features & 10 questionnaire items (Q1--Q10) \\
Feature encoding & Centered Likert encoding ($[-1,1]$) \\
Train--test split & 70\% / 30\% \\
Sampling strategy & Stratified by gender \\
Randomization & Fixed seed with repeated splits \\
Classifier & Linear SVC \\
Regularization & $C$ via cross-validation \\
Dimensionality reduction & PCA (analysis + training) \\
Evaluation & Mean performance over runs \\
Implementation & Python 3.10, \texttt{scikit-learn} \\
\bottomrule
\end{tabular}
\end{table}

For each split, the SVC is trained using only the training subset. The regularization parameter $C$ is selected through cross-validation on the training data, preventing information leakage from the test set. The final model is then retrained on the full training subset using the selected parameter and evaluated on the held-out test data.

To assess the stability of model performance, the training and evaluation procedure is repeated across multiple randomized stratified splits. Performance metrics, including classification accuracy and ROC--AUC, are computed for each run and aggregated to obtain mean values and variability estimates. ROC--AUC denotes the area under the receiver operating characteristic curve and measures the ability of the classifier to distinguish between the two classes across decision thresholds \cite{}. Confidence intervals for classification accuracy are estimated using both pooled binomial Wilson intervals and empirical intervals derived from repeated splits.

The regenerated dataset is used only to stabilize model training and support repeated evaluation under limited sample size and class imbalance. All conclusions are based on aggregated results across repeated evaluations rather than on individual regenerated observations. All experiments are implemented in Python~3.10 using the \texttt{scikit-learn} library.
\section{Results}
\label{sec:results}

\subsection{Overall Performance}

\begin{figure}[!ht]
\centering
\includegraphics[width=\linewidth]{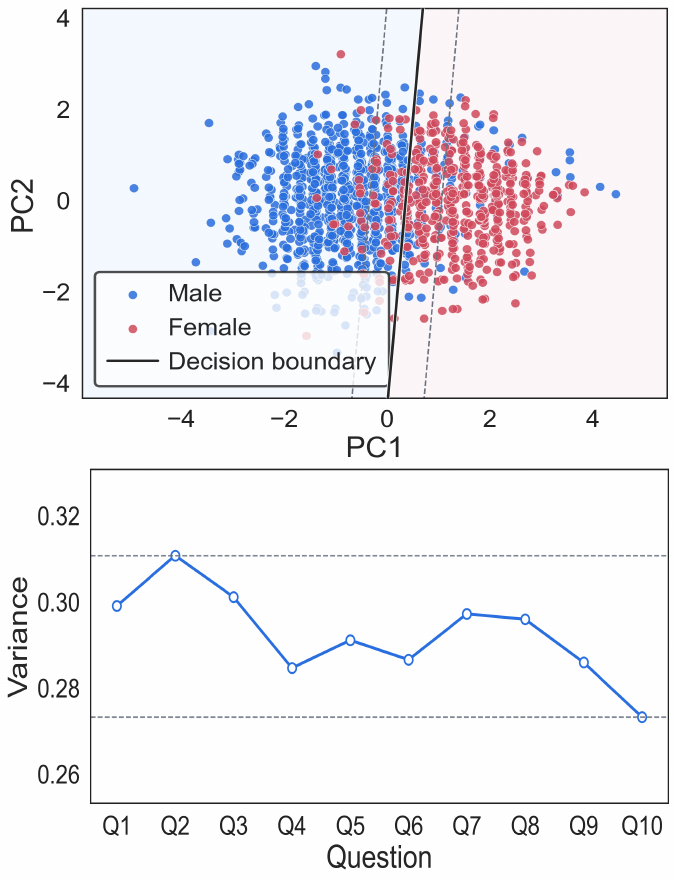}
\caption{Linear SVM classification of gender from survey responses. 
\textbf{Top:} Decision boundary in PCA space showing structured but overlapping separation. 
\textbf{Bottom:} Coefficient magnitudes indicating the relative importance of questionnaire items.}
\label{fig:fig3}
\end{figure}

We evaluate the SVC on the task of distinguishing male and female response patterns using the encoded questionnaire data. Figure~\ref{fig:fig3} shows the classification results. The top panel presents the learned decision boundary projected onto the first two principal components, illustrating a structured but overlapping separation between the groups. The bottom panel shows the magnitude of the model coefficients, providing an interpretable ranking of the questionnaire items.

When trained in the original 10-dimensional feature space (Q1--Q10), the classifier achieves a test accuracy of $85.24\%$ and a ROC--AUC of $0.9205$. The mean absolute decision margin is $1.72$, indicating a clear separation between the groups. This separation arises from consistent differences in multivariate response patterns rather than any single questionnaire item, suggesting that group-specific perceptions are reflected in structured combinations of responses. The overall performance metrics and PCA variance statistics are reported in Table~\ref{tab:overall_perf}.

\begin{table}[!ht]
\centering
\caption{Overall classification performance and PCA variance statistics.}
\label{tab:overall_perf}
\renewcommand{\arraystretch}{1.05}
\begin{tabular}{lc}
\hline
\textbf{Metric} & \textbf{Value} \\
\hline
Test accuracy & $85.24\%$ \\
ROC--AUC & $0.9205$ \\
Mean $|\mathrm{margin}|$ & $1.72$ \\
Std.\ $|\mathrm{margin}|$ & $1.13$ \\
PCA variance (PC1) & $20.44\%$ \\
PCA variance (PC2) & $9.75\%$ \\
Total PCA variance & $30.19\%$ \\
\hline
\end{tabular}
\end{table}

To examine the role of dominant variance structure, responses are projected onto the first two principal components, which together explain $30.19\%$ of the total variance. Despite this reduction, classification performance remains comparable to that in the full feature space, suggesting that the main variance directions capture a substantial portion of the discriminative information. This indicates that the observed separation reflects stable patterns in the data rather than random variation.

\subsection{Performance Stability}

To assess robustness to data partitioning, the training and evaluation procedure is repeated across $R=20$ randomized stratified train--test splits ($70\%/30\%$). The classifier achieves a mean accuracy of $86.29\%$ with a standard deviation of $0.63\%$, and an empirical $95\%$ interval of $[85.21\%,\,87.26\%]$, indicating consistent performance across different splits. This stability suggests that the separation is not sensitive to a particular data partition, but reflects a reliable pattern in the data. The repeated-split stability statistics are summarized in Table~\ref{tab:stat_comparison}.

\begin{table}[!h]
\centering
\caption{Classification stability over repeated splits.}
\label{tab:stat_comparison}
\small
\renewcommand{\arraystretch}{1.05}
\begin{tabular}{@{}lc@{}}
\toprule
\textbf{Statistic} & \textbf{Value} \\
\midrule
Mean acc.        & 86.29\% \\
Std.\ dev.       & 0.63\% \\
Min acc.         & 85.11\% \\
Max acc.         & 87.51\% \\
Empirical 95\%   & [85.21\%, 87.26\%] \\
\bottomrule
\end{tabular}
\end{table}

\begin{figure}[!ht]
\centering
\includegraphics[width=\linewidth]{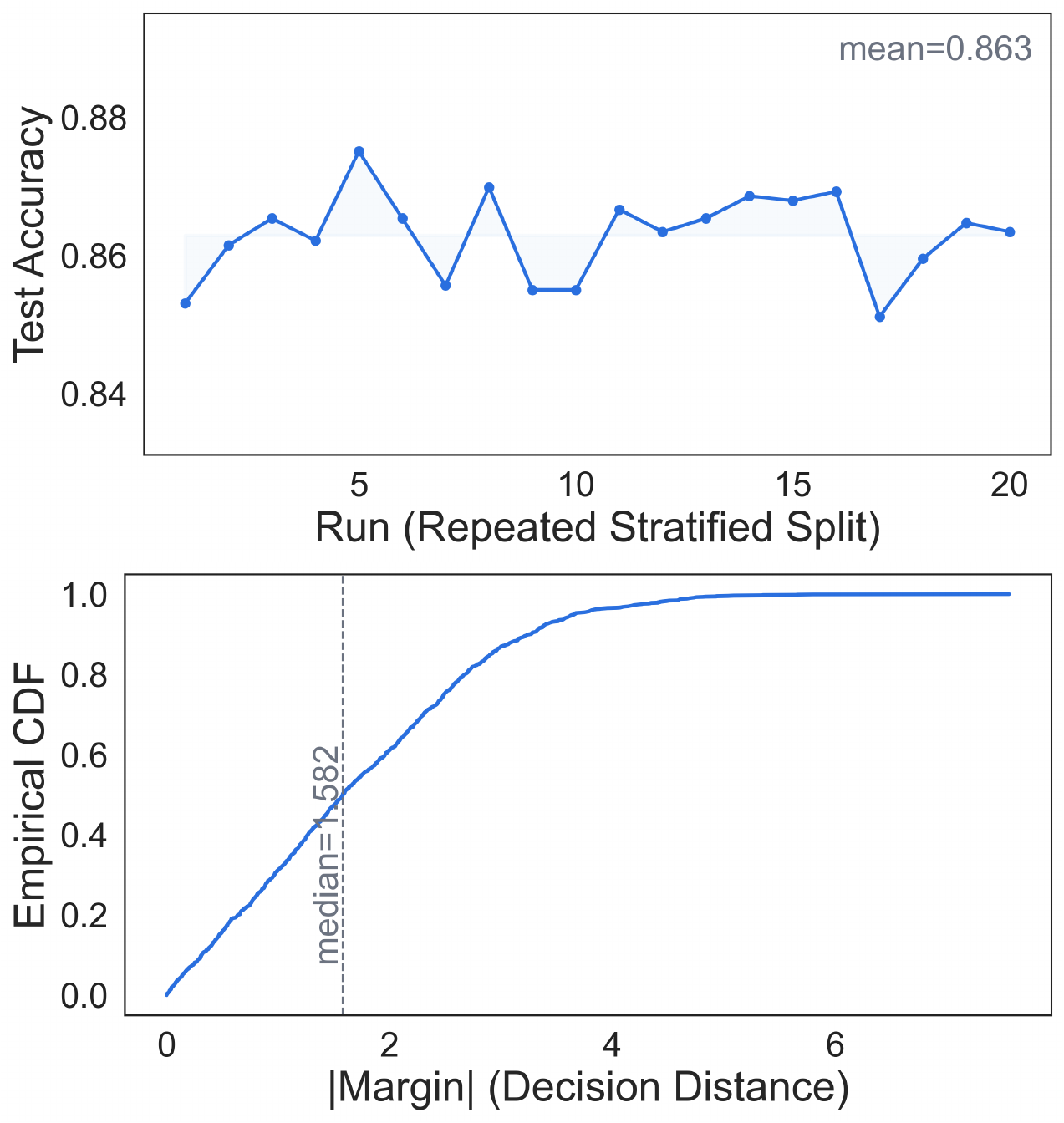}
\caption{Stability of the linear SVM classifier. 
\textbf{Top:} Test accuracy across repeated stratified splits. 
\textbf{Bottom:} Empirical distribution of decision margins.}
\label{fig:fig4}
\end{figure}

Figure~\ref{fig:fig4} further illustrates this behavior. The top panel shows classification accuracy across repeated splits, with only minor variation around the mean. The bottom panel shows the empirical distribution of decision margins, where most samples lie at a non-zero distance from the decision boundary, indicating confident and stable classification.

Class-wise performance for a representative test split is summarized in Table~\ref{tab:class_metrics}. Recall values are comparable across classes, indicating that the classifier does not systematically under-identify either group.

\begin{table}[!h]
\centering
\caption{Class-wise precision, recall, and F1-score for a representative split.}
\label{tab:class_metrics}
\small
\renewcommand{\arraystretch}{1.05}
\begin{tabular}{@{}lccc@{}}
\toprule
\textbf{Class} & \textbf{Precision} & \textbf{Recall} & \textbf{F1-score} \\
\midrule
Male   & 0.92 & 0.85 & 0.88 \\
Female & 0.74 & 0.86 & 0.80 \\
\bottomrule
\end{tabular}
\end{table}

\subsection{Feature Contributions}

The magnitude of the linear SVM coefficients provides an interpretable ranking of questionnaire items contributing to the classification rule. The resulting coefficient-based ranking is reported in Table~\ref{tab:feature_rank}. In this table, $w_j$ denotes the learned coefficient associated with questionnaire item $j$, and $|w_j|$ denotes its absolute magnitude. Larger values of $|w_j|$ indicate stronger contribution to the linear decision rule.

\begin{table}[!h]
\centering
\caption{Ranking of questionnaire items based on the magnitude of linear SVM coefficients.}
\label{tab:feature_rank}
\renewcommand{\arraystretch}{1.05}
\begin{tabular}{lcc}
\hline
\textbf{Question} & \textbf{Weight $w_j$} & \textbf{$|w_j|$} \\
\hline
Q1  & $0.252$ & $0.252$ \\
Q2  & $0.239$ & $0.239$ \\
Q8  & $0.237$ & $0.237$ \\
Q3  & $0.232$ & $0.232$ \\
Q9  & $0.231$ & $0.231$ \\
Q7  & $0.230$ & $0.230$ \\
Q5  & $0.183$ & $0.183$ \\
Q6  & $0.170$ & $0.170$ \\
Q4  & $0.166$ & $0.166$ \\
Q10 & $0.140$ & $0.140$ \\
\hline
\end{tabular}
\end{table}

The highest-weighted items, Q1 and Q2, relate to students’ engagement and overall experience with the flipped classroom. Items associated with instructional organization and facilitation, such as Q8 and Q9, also contribute strongly. This pattern suggests that the observed separation is driven by a combination of engagement-related and instructional factors, rather than by any single questionnaire item.

\subsection{Comparison with Alternative Classifiers}

To verify that the observed separation is not model-specific, the linear SVM is compared with logistic regression, random forest, and gradient-boosted trees (XGBoost). All models are trained using the same preprocessing pipeline and repeated stratified train--test splits. The comparison across alternative classifiers is summarized in Table~\ref{tab:model_comparison}.

\begin{table}[!h]
\centering
\caption{Classification accuracy across classifiers, reported as mean accuracy over repeated splits.}
\label{tab:model_comparison}
\small
\renewcommand{\arraystretch}{1.05}
\begin{tabular}{@{}lcc@{}}
\toprule
\textbf{Classifier} & \textbf{Acc. (\%)} & \textbf{95\% CI} \\
\midrule
Logistic Regression & 86.22 & [85.0, 87.2] \\
Random Forest       & 86.27 & [84.7, 87.4] \\
XGBoost             & 85.89 & [84.4, 87.3] \\
\textbf{Linear SVM} & \textbf{86.29} & \textbf{[85.21, 87.26]} \\
\bottomrule
\end{tabular}
\end{table}

All classifiers achieve comparable performance with overlapping confidence intervals. This consistency indicates that the classification pattern reflects a stable structure in the survey data rather than dependence on a specific model. The similar performance of linear and non-linear models further suggests that these differences are organized along approximately linear directions in the questionnaire feature space.

\subsection{Limitations and Scope of Interpretation}

Survey responses reflect students’ perceptions and self-reported experiences of the flipped classroom, rather than direct measures of learning outcomes such as examination performance. As a result, the classification results should be interpreted as differences in reported perceptions, not as differences in academic achievement.

Questionnaire-based data may also be subject to response bias, limited expressiveness, or incomplete information. If responses across groups were largely similar or dominated by noise, the resulting feature space would exhibit substantial overlap and classification performance would approach random. The results therefore depend on the presence of sufficiently informative and structured response patterns in the survey data.

Although the use of an interpretable model allows group differences to be linked to specific questionnaire items, the findings remain dependent on the quality and design of the survey instrument. Accordingly, the results should be understood as patterns in reported perceptions rather than definitive conclusions about underlying learning behavior.

The formulation of the problem as a gender-based classification task should also be interpreted with caution. The objective is to identify structured differences in response patterns, not to imply inherent or causal distinctions between groups. These differences may reflect a combination of pedagogical, contextual, and response-related factors, and the present analysis does not establish their underlying causes.

Finally, while the analysis identifies key dimensions such as engagement and instructional clarity, it does not directly establish how these factors translate into specific pedagogical interventions or learning outcomes. Further work is required to connect these patterns to actionable instructional design and to evaluate their relationship with student achievement.

\section{Conclusion}
\label{sec:conclusion}

This study analyzed student perceptions of a flipped introductory linear algebra course using an interpretable machine-learning framework. The results show that multivariate survey responses contain consistent and stable structure across gender groups, which can be captured using a simple linear model. The robustness of this pattern across repeated evaluations suggests that it reflects meaningful variation in the data rather than artifacts of sampling or model choice.

The identified differences are distributed across multiple pedagogical dimensions, particularly those related to engagement and instructional organization. This indicates that variation in student perception arises from combinations of factors rather than isolated questionnaire responses. These findings suggest that while flipped classroom pedagogy is broadly effective, targeted refinements may better support diverse student experiences.

Because this analysis is based on perception-driven survey data, the results should not be interpreted as direct evidence of differences in academic performance. Future work may examine how these perception patterns relate to learning outcomes, classroom participation, and longer-term student achievement. Overall, this study highlights the value of interpretable machine learning for understanding structured patterns in educational data and for supporting more inclusive instructional design.

\section*{Conflict of Interest Statement}

The authors declare that they have no conflict of interest.
{
    \small
    \bibliographystyle{ieeenat_fullname}
    \bibliography{main}

\begin{thebibliography}{34}
\providecommand{\natexlab}[1]{#1}
\providecommand{\url}[1]{\texttt{#1}}
\expandafter\ifx\csname urlstyle\endcsname\relax
  \providecommand{\doi}[1]{doi: #1}\else
  \providecommand{\doi}{doi: \begingroup \urlstyle{rm}\Url}\fi

\bibitem[Bergmann and Sams(2023)]{bergmann2023flip}
Jon Bergmann and Aaron Sams.
\newblock \emph{Flip Your Classroom, Revised Edition: Reach Every Student in Every Class Every Day}.
\newblock ASCD, 2023.

\bibitem[Bishop and Verleger(2013)]{bishop2013flipped}
Jacob Bishop and Matthew~A. Verleger.
\newblock The flipped classroom: A survey of the research.
\newblock In \emph{2013 ASEE Annual Conference \& Exposition}, pages 23--1200, 2013.

\bibitem[Carifio and Perla(2008)]{carifio2008resolving}
James Carifio and Rocco Perla.
\newblock Resolving the 50-year debate around using and misusing {Likert} scales.
\newblock \emph{Medical Education}, 42\penalty0 (12):\penalty0 1150--1152, 2008.

\bibitem[Cortes and Vapnik(1995)]{cortes1995support}
Corinna Cortes and Vladimir Vapnik.
\newblock Support-vector networks.
\newblock \emph{Machine Learning}, 20:\penalty0 273--297, 1995.

\bibitem[Eddy and Hogan(2014)]{eddy2014getting}
Sarah~L. Eddy and Kelly~A. Hogan.
\newblock Getting under the hood: How and for whom does increasing course structure work?
\newblock \emph{CBE---Life Sciences Education}, 13\penalty0 (3):\penalty0 453--468, 2014.

\bibitem[Flood et~al.(2025)Flood, England, and Grawemeyer]{flood2025predicting}
Daniel Flood, Matthew England, and Beate Grawemeyer.
\newblock Predicting at-risk programming students in small imbalanced datasets using synthetic data.
\newblock In \emph{International Conference on Artificial Intelligence in Education}, pages 427--432. Springer, 2025.

\bibitem[Freeman et~al.(2014)Freeman, Eddy, McDonough, Smith, Okoroafor, Jordt, and Wenderoth]{freeman2014active}
Scott Freeman, Sarah~L. Eddy, Miles McDonough, Michelle~K. Smith, Nnadozie Okoroafor, Hannah Jordt, and Mary~Pat Wenderoth.
\newblock Active learning increases student performance in science, engineering, and mathematics.
\newblock \emph{Proceedings of the National Academy of Sciences}, 111\penalty0 (23):\penalty0 8410--8415, 2014.

\bibitem[Hardebolle et~al.(2022)Hardebolle, Verma, Tormey, and Deparis]{hardebolle2022gender}
C{\'e}cile Hardebolle, Himanshu Verma, Roland Tormey, and Simone Deparis.
\newblock Gender, prior knowledge, and the impact of a flipped linear algebra course for engineers over multiple years.
\newblock \emph{Journal of Engineering Education}, 111\penalty0 (3):\penalty0 554--574, 2022.

\bibitem[He et~al.(2016)He, Holton, Farkas, and Warschauer]{he2016effects}
Wenliang He, Alycia Holton, George Farkas, and Mark Warschauer.
\newblock The effects of flipped instruction on out-of-class study time, exam performance, and student perceptions.
\newblock \emph{Learning and Instruction}, 45:\penalty0 61--71, 2016.

\bibitem[Jolliffe(2011)]{jolliffe2011principal}
Ian Jolliffe.
\newblock Principal component analysis.
\newblock In \emph{International Encyclopedia of Statistical Science}, pages 1094--1096. Springer, 2011.

\bibitem[Karjanto and Acelajado(2022)]{karjanto2022sustainable}
Natanael Karjanto and Maxima~J. Acelajado.
\newblock Sustainable learning, cognitive gains, and improved attitudes in {College Algebra} flipped classrooms.
\newblock \emph{Sustainability}, 14\penalty0 (19):\penalty0 12500, 2022.

\bibitem[Karjanto and Simon(2019)]{karjanto2019english}
Natanael Karjanto and Lois Simon.
\newblock English-medium instruction {Calculus} in {Confucian-Heritage Culture}: Flipping the class or overriding the culture?
\newblock \emph{Studies in Educational Evaluation}, 63:\penalty0 122--135, 2019.

\bibitem[Kohavi(1995)]{kohavi1995study}
Ron Kohavi.
\newblock A study of cross-validation and bootstrap for accuracy estimation and model selection.
\newblock In \emph{Proceedings of the 14th International Joint Conference on Artificial Intelligence}, pages 1137--1145, 1995.

\bibitem[Lage et~al.(2000)Lage, Platt, and Treglia]{lage2000inverting}
Maureen~J. Lage, Glenn~J. Platt, and Michael Treglia.
\newblock Inverting the classroom: A gateway to creating an inclusive learning environment.
\newblock \emph{The Journal of Economic Education}, 31\penalty0 (1):\penalty0 30--43, 2000.

\bibitem[Livieris et~al.(2023)Livieris, Karacapilidis, Domalis, and Tsakalidis]{livieris2023advanced}
Ioannis~E. Livieris, Nikos Karacapilidis, Georgios Domalis, and Dimitris Tsakalidis.
\newblock An advanced explainable and interpretable {ML}-based framework for educational data mining.
\newblock In \emph{International Conference in Methodologies and Intelligent Systems for Technology Enhanced Learning}, pages 87--96. Springer, 2023.

\bibitem[Love et~al.(2014)Love, Hodge, Grandgenett, and Swift]{love2014student}
Betty Love, Angie Hodge, Neal Grandgenett, and Andrew~W. Swift.
\newblock Student learning and perceptions in a flipped linear algebra course.
\newblock \emph{International Journal of Mathematical Education in Science and Technology}, 45\penalty0 (3):\penalty0 317--324, 2014.

\bibitem[Murphy et~al.(2016)Murphy, Chang, and Suaray]{murphy2016student}
Julia Murphy, Jen-Mei Chang, and Kagba Suaray.
\newblock Student performance and attitudes in a collaborative and flipped linear algebra course.
\newblock \emph{International Journal of Mathematical Education in Science and Technology}, 47\penalty0 (5):\penalty0 653--673, 2016.

\bibitem[Nasir et~al.(2020)Nasir, Alaudin, Ismail, Ali, Faudzi, Yusuff, and Pozi]{nasir2020effectiveness}
Mohd Azrin~Mohd Nasir, Ros~Idayuwati Alaudin, Suzila Ismail, Nor A'tikah~Mat Ali, Farah Nadia~Mohd Faudzi, Noraini Yusuff, and Muhammad Syafiq~Mohd Pozi.
\newblock The effectiveness of flipped classroom strategy on self-directed learning among undergraduate mathematics students.
\newblock \emph{Practitioner Research}, 2:\penalty0 61--81, 2020.

\bibitem[Noble(2006)]{noble2006support}
William~S. Noble.
\newblock What is a support vector machine?
\newblock \emph{Nature Biotechnology}, 24\penalty0 (12):\penalty0 1565--1567, 2006.

\bibitem[Norman(2010)]{norman2010likert}
Geoff Norman.
\newblock {Likert} scales, levels of measurement and the laws of statistics.
\newblock \emph{Advances in Health Sciences Education}, 15\penalty0 (5):\penalty0 625--632, 2010.

\bibitem[Novak et~al.(2017)Novak, Kensington-Miller, and Evans]{novak2017flip}
Julia Novak, Barbara Kensington-Miller, and Tanya Evans.
\newblock Flip or flop? {Students'} perspectives of a flipped lecture in mathematics.
\newblock \emph{International Journal of Mathematical Education in Science and Technology}, 48\penalty0 (5):\penalty0 647--658, 2017.

\bibitem[O'Flaherty and Phillips(2015)]{o2015use}
Jacqueline O'Flaherty and Craig Phillips.
\newblock The use of flipped classrooms in higher education: A scoping review.
\newblock \emph{The Internet and Higher Education}, 25:\penalty0 85--95, 2015.

\bibitem[Pedregosa et~al.(2011)Pedregosa, Varoquaux, Gramfort, Michel, Thirion, Grisel, Blondel, Prettenhofer, Weiss, Dubourg, Vanderplas, Passos, Cournapeau, Brucher, Perrot, and Duchesnay]{pedregosa2011scikit}
Fabian Pedregosa, Ga{\"e}l Varoquaux, Alexandre Gramfort, Vincent Michel, Bertrand Thirion, Olivier Grisel, Mathieu Blondel, Peter Prettenhofer, Ron Weiss, Vincent Dubourg, Jake Vanderplas, Alexandre Passos, David Cournapeau, Matthieu Brucher, Matthieu Perrot, and {\'E}douard Duchesnay.
\newblock Scikit-learn: Machine learning in {Python}.
\newblock \emph{Journal of Machine Learning Research}, 12:\penalty0 2825--2830, 2011.

\bibitem[Prince(2004)]{prince2004does}
Michael Prince.
\newblock Does active learning work? {A} review of the research.
\newblock \emph{Journal of Engineering Education}, 93\penalty0 (3):\penalty0 223--231, 2004.

\bibitem[Ranjeeth et~al.(2020)Ranjeeth, Latchoumi, and Paul]{ranjeeth2020survey}
Sama Ranjeeth, Thamarai~Pugazhendhi Latchoumi, and P.~Victer Paul.
\newblock A survey on predictive models of learning analytics.
\newblock \emph{Procedia Computer Science}, 167:\penalty0 37--46, 2020.

\bibitem[Romero and Ventura(2020)]{romero2020educational}
Cristobal Romero and Sebastian Ventura.
\newblock Educational data mining and learning analytics: An updated survey.
\newblock \emph{Wiley Interdisciplinary Reviews: Data Mining and Knowledge Discovery}, 10\penalty0 (3):\penalty0 e1355, 2020.

\bibitem[Santana-Perera et~al.(2025)Santana-Perera, Garc{\'\i}a-Barcel{\'o}, Gonz{\'a}lez~Arcas, and Gil]{santana2025exploring}
Beatriz Santana-Perera, Carmen Garc{\'\i}a-Barcel{\'o}, Mauricio Gonz{\'a}lez~Arcas, and David Gil.
\newblock Exploring predictive insights on student success using explainable machine learning: A synthetic data study.
\newblock \emph{Information}, 16\penalty0 (9):\penalty0 763, 2025.

\bibitem[Strelan et~al.(2020{\natexlab{a}})Strelan, Osborn, and Palmer]{strelan2020flipped}
Peter Strelan, Amanda Osborn, and Edward Palmer.
\newblock The flipped classroom: A meta-analysis of effects on student performance across disciplines and education levels.
\newblock \emph{Educational Research Review}, 30:\penalty0 100314, 2020{\natexlab{a}}.

\bibitem[Strelan et~al.(2020{\natexlab{b}})Strelan, Osborn, and Palmer]{strelan2020student}
Peter Strelan, Amanda Osborn, and Edward Palmer.
\newblock Student satisfaction with courses and instructors in a flipped classroom: A meta-analysis.
\newblock \emph{Journal of Computer Assisted Learning}, 36\penalty0 (3):\penalty0 295--314, 2020{\natexlab{b}}.

\bibitem[Suthaharan(2016)]{suthaharan2016support}
Shan Suthaharan.
\newblock Support vector machine.
\newblock In \emph{Machine Learning Models and Algorithms for Big Data Classification: Thinking with Examples for Effective Learning}, pages 207--235. Springer, Boston, MA, USA, 2016.

\bibitem[Takrouni et~al.(2024)Takrouni, Neji, and Jguirim]{takrouni2024data}
Manel Takrouni, Wissal Neji, and Nabil Jguirim.
\newblock A data-driven flipped learning model for enhanced student outcomes.
\newblock In \emph{Proceedings of the 20th International CDIO Conference}, 2024.
\newblock Hosted by Ecole Supérieure Privée d'Ingénierie et de Technologies (ESPRIT), Tunis, Tunisia, June 10--13, 2024.

\bibitem[Talbert(2014)]{talbert2014inverting}
Robert Talbert.
\newblock Inverting the linear algebra classroom.
\newblock \emph{PRIMUS}, 24\penalty0 (5):\penalty0 361--374, 2014.

\bibitem[Whitelock-Wainwright et~al.(2019)Whitelock-Wainwright, Ga{\v{s}}evi{\'c}, Tejeiro, Tsai, and Bennett]{whitelock2019student}
Alexander Whitelock-Wainwright, Dragan Ga{\v{s}}evi{\'c}, Ricardo Tejeiro, Yi-Shan Tsai, and Kate Bennett.
\newblock The student expectations of learning analytics questionnaire.
\newblock \emph{Journal of Computer Assisted Learning}, 35\penalty0 (5):\penalty0 633--666, 2019.

\bibitem[Zhao et~al.(2025)Zhao, Yuan, Luo, Xie, Zhang, Quan, Yuan, Wang, and Zhang]{zhao2025large}
Jianpeng Zhao, Chenyu Yuan, Weiming Luo, Haoling Xie, Guangwei Zhang, Steven~Jige Quan, Zixuan Yuan, Pengyang Wang, and Denghui Zhang.
\newblock Large language models as virtual survey respondents: Evaluating sociodemographic response generation.
\newblock \emph{arXiv preprint arXiv:2509.06337}, 2025.

\end{thebibliography}
}

\clearpage
\setcounter{page}{1}
\maketitlesupplementary

\section*{Appendix: Survey Instrument for Flipped Learning in Introductory Linear Algebra}
\label{sec:supp}

This appendix presents the questionnaire used to collect students’ perceptions and attitudes toward flipped learning in an introductory linear algebra course. The survey was administered at the end of the semester.


\subsection*{Student Information}

Please respond to the following items by marking the appropriate boxes.

\medskip

\noindent
\begin{tabularx}{\linewidth}{@{}p{0.24\linewidth}X@{}}
\textbf{Gender:} & $\Box$ Male \quad $\Box$ Female \\[0.4em]

\textbf{Study year:} & 
$\Box$ Freshman (Year 1) \qquad \newline
$\Box$ Sophomore (Year 2) \quad \newline
$\Box$ Junior (Year 3) \quad \newline
$\Box$ Senior (Year 4) \quad \newline
$\Box$ Other \vspace*{0.1cm} \\[0.4em] 

\textbf{School:} & 
$\Box$ \!Natural Science \quad
$\Box$ Engineering \qquad 
$\Box$ \;Computer Science \quad \,
$\Box$ \; Other \vspace*{0.1cm} \\[0.4em]

\textbf{Expected grade:} & 
$\Box$ A+ \qquad
$\Box$ A \qquad
$\Box$ B+ \qquad
$\Box$ B \qquad 
$\Box$ C+ \qquad 
$\Box$ C \qquad
$\Box$ D+ \qquad
$\Box$ D \qquad
$\Box$ \, F \qquad \vspace*{0.1cm} \\[0.4em]

\textbf{Current GPA:} & 
\underline{\hspace{3cm}} \hspace{4cm}
{\small Please estimate if you do not know; \newline provide a single numeric value, e.g., 3.78.}
\end{tabularx}

\newpage
\subsection*{Flipped Learning Perception Questionnaire}

The following statements are designed to collect your opinions and attitudes toward flipped learning in Linear Algebra. Each item has five possible responses ranging from \emph{Strongly disagree} to \emph{Strongly agree}. 
If you have no opinion, please select \emph{Neutral}. Please read each statement carefully and select the response that best reflects your level of agreement.

\medskip

\begin{center}
\small
\setlength{\tabcolsep}{2.5pt}
\renewcommand{\arraystretch}{1.1}
\begin{tabularx}{\linewidth}{>{\raggedright\arraybackslash}X *{5}{c}}
\hline
\textbf{Statement} & \textbf{SD} & \textbf{D} & \textbf{N} & \textbf{A} & \textbf{SA} \\
\hline
I prefer flipped classroom more than the traditional classroom. 
& $\square$ & $\square$ & $\square$ & $\square$ & $\square$ \\

I enjoy flipped classroom more than the traditional classroom. 
& $\square$ & $\square$ & $\square$ & $\square$ & $\square$ \\

I watch recorded video lectures before coming to the class. 
& $\square$ & $\square$ & $\square$ & $\square$ & $\square$ \\

I like to participate in problem-solving activities during the class. 
& $\square$ & $\square$ & $\square$ & $\square$ & $\square$ \\

I understand better the course material when it is delivered in the flipped learning format. 
& $\square$ & $\square$ & $\square$ & $\square$ & $\square$ \\

My instructor provided clear instructions on how to participate in course activities. 
& $\square$ & $\square$ & $\square$ & $\square$ & $\square$ \\

My instructor successfully implemented flipped learning for this course. 
& $\square$ & $\square$ & $\square$ & $\square$ & $\square$ \\

My instructor successfully encouraged the students to participate actively during the class. 
& $\square$ & $\square$ & $\square$ & $\square$ & $\square$ \\

My instructor organized learning activities for effective and efficient use of time. 
& $\square$ & $\square$ & $\square$ & $\square$ & $\square$ \\

My instructor is helpful in responding to my questions or problems. 
& $\square$ & $\square$ & $\square$ & $\square$ & $\square$ \\
\hline
\end{tabularx}
\end{center}

\medskip

\noindent\textbf{Additional comments (optional):}  
\vspace{1cm}

\end{document}